

\magnification=\magstep1
\baselineskip =5mm
\lineskiplimit =1.0mm
\lineskip =1.0mm

\long\def\comment#1{}

\long\def\blankout #1\eb{}
\def\noblankout{\def\blankout{}\def\eb{}}

\let\properlbrack=\lbrack
\let\properrbrack=\rbrack
\def\ordcomma{,}
\def\ordcolon{:}
\def\ordsemicolon{;}
\def\ordleftparen{(}
\def\ordrightparen{)}
\def\ordleftbrack{\properlbrack}
\def\ordrightbrack{\properrbrack}
\def\rmcomma{\ifmmode ,\else \/{\rm ,}\fi}
\def\rmcolon{\ifmmode :\else \/{\rm :}\fi}
\def\rmsemicolon{\ifmmode ;\else \/{\rm ;}\fi}
\def\rmleftparen{\ifmmode (\else \/{\rm (}\fi}
\def\rmrightparen{\ifmmode )\else \/{\rm )}\fi}
\def\rmleftbrack{\ifmmode \properlbrack\else \/{\rm \properlbrack}\fi}
\def\rmrightbrack{\ifmmode \properrbrack\else \/{\rm \properrbrack}\fi}
\catcode`,=\active 
\catcode`:=\active 
\catcode`;=\active 
\catcode`(=\active 
\catcode`)=\active 
\catcode`[=\active 
\catcode`]=\active 
\let,=\ordcomma
\let:=\ordcolon
\let;=\ordsemicolon
\let(=\ordleftparen
\let)=\ordrightparen
\let[=\ordleftbrack
\let]=\ordrightbrack
\let\lbrack=\ordleftbrack
\let\rbrack=\ordrightbrack
\def\rmpunctuation{
\let,=\rmcomma
\let:=\rmcolon
\let;=\rmsemicolon
\let(=\rmleftparen
\let)=\rmrightparen
\let[=\rmleftbrack
\let]=\rmrightbrack
\let\lbrack=\rmleftbrack
\let\rbrack=\rmrightbrack}

\def\writemonth#1{\ifcase#1
\or January\or February\or March\or April\or May\or June\or July%
\or August\or September\or October\or November\or December\fi}

\newcount\mins
\newcount\minmodhour
\newcount\hour
\newcount\hourinmin
\newcount\ampm
\newcount\ampminhour
\newcount\hourmodampm
\def\writetime#1{%
\mins=#1%
\hour=\mins \divide\hour by 60
\hourinmin=\hour \multiply\hourinmin by -60
\minmodhour=\mins \advance\minmodhour by \hourinmin
\ampm=\hour \divide\ampm by 12
\ampminhour=\ampm \multiply\ampminhour by -12
\hourmodampm=\hour \advance\hourmodampm by \ampminhour
\ifnum\hourmodampm=0 12\else \number\hourmodampm\fi
:\ifnum\minmodhour<10 0\number\minmodhour\else \number\minmodhour\fi
\ifodd\ampm p.m.\else a.m.\fi
}

\font\tenrm=cmr10
\font\smallcaps=cmcsc10
\font\eightrm=cmr8
\font\ninerm=cmr9
\font\sixrm=cmr6
\font\eightbf=cmbx8
\font\sixbf=cmbx6
\font\eightit=cmti8
\font\eightsl=cmsl8
\font\eighti=cmmi8
\font\eightsy=cmsy8
\font\eightex=cmex10 at 8pt
\font\sixi=cmmi6
\font\sixsy=cmsy6
\font\ninesy=cmsy9
\font\seventeenrm=cmr17
 1
\font\twelverm=cmr10 scaled \magstep2
\font\seventeeni=cmmi10 scaled \magstep3
\font\twelvei=cmmi10 scaled \magstep2
\font\seventeensy=cmsy10 scaled \magstep3
\font\twelvesy=cmsy10 at 12pt
\font\seventeenex=cmex10 scaled \magstep3
\font\seventeenbf=cmbx10 scaled \magstep3
\catcode`@=11
\def\eightbig#1{{\hbox{$\textfont0=\ninerm\textfont2=\ninesy
\left#1\vbox to6.5pt{}\right.\n@space$}}}
\catcode`@=12
\def\eightpoint{\eightrm \normalbaselineskip=4.5 mm%
\textfont0=\eightrm \scriptfont0=\sixrm \scriptscriptfont0=\fiverm%
\def\rm{\fam0 \eightrm}%
\textfont1=\eighti \scriptfont1=\sixi \scriptscriptfont1=\fivei%
\def\mit{\fam1 } \def\oldstyle{\fam1 \eighti}%
\textfont2=\eightsy \scriptfont2=\sixsy \scriptscriptfont2=\fivesy%
\def\cal{\fam2 }%
\textfont3=\eightex \scriptfont3=\eightex \scriptscriptfont3=\eightex%
\def\bf{\fam\bffam\eightbf} \textfont\bffam\eightbf
\scriptfont\bffam=\sixbf \scriptscriptfont\bffam=\fivebf
\def\it{\fam\itfam\eightit} \textfont\itfam\eightit
\def\sl{\fam\slfam\eightsl} \textfont\slfam\eightsl
\let\big=\eightbig \normalbaselines\rm
\def\caps##1{\nottencaps{##1}}
}
\def\seventeenpoint{\seventeenrm \baselineskip=5.5mm%
\textfont0=\seventeenrm \scriptfont0=\twelverm \scriptscriptfont0=\sevenrm%
\def\rm{\fam0 \seventeenrm}%
\textfont1=\seventeeni \scriptfont1=\twelvei \scriptscriptfont1=\seveni%
\def\mit{\fam1 } \def\oldstyle{\fam1 \seventeeni}%
\textfont2=\seventeensy \scriptfont2=\twelvesy \scriptscriptfont2=\sevensy%
\def\cal{\fam2 }%
\textfont3=\seventeenex \scriptfont3=\seventeenex%
\scriptscriptfont3=\seventeenex%
\def\bf{\fam\bffam\seventeenbf} \textfont\bffam\seventeenbf
}

\def\setheadline #1\\ #2 \par{\headline={\ifnum\pageno=1 
\hfil
\else \eightpoint \noindent
\ifodd\pageno \hfil \caps{#2}\hfil \else
\hfil \caps{#1}\hfil \fi\fi}}

\def\beginsection{} 
\def\datedversion{\footline={\ifnum\pageno=1 \fiverm \hfil
Typeset using plain-\TeX\ on
\writemonth\month\ \number\day, \number\year\ at \writetime{\time}\hfil 
\else \tenrm \hfil \folio \hfil \fi}
\def\tempsetheadline##1{\headline={\ifnum\pageno=1 
\hfil
\else \eightpoint \noindent
\writemonth\month\ \number\day, \number\year,
\writetime{\time}\hfil ##1\fi}}
\def\firstbeginsection##1\par{\bigskip\vskip\parskip
\message{##1}\centerline{\caps{##1}}\nobreak\smallskip\noindent
\tempsetheadline{##1}} \def\beginsection##1\par{\vskip0pt
plus.3\vsize\penalty-250 \vskip0pt plus-.3\vsize\bigskip\vskip\parskip
\message{##1}\centerline{\caps{##1}}\nobreak\smallskip\noindent
\tempsetheadline{##1}}}

\def\finalversion{\footline={\ifnum\pageno=1 \eightrm \hfil 
This paper is in final form.\hfil 
\else \tenrm \hfil \folio \hfil \fi}}

\def\preliminaryversion{\footline={\ifnum\pageno=1 \eightrm \hfil 
Preliminary Version.\hfil 
\else \tenrm \hfil \folio \hfil \fi}}

\def\moreproclaim{\par}
\def\Head #1: {\medskip\noindent{\it #1}:\enspace}
\def\Proof: {\Head Proof: }
\def\Proofof #1: {\Head Proof of #1: }
\def\endproof{\nobreak\hfill$\sqr$\bigskip\goodbreak}
\def\itemi{\item{i)}}
\def\itemii{\item{ii)}}
\def\itemiii{\item{iii)}}

\def\ds#1{{\displaystyle{#1}}}
\def\ts#1{{\textstyle{#1}}}

\def\Abstract\par#1\par{\centerline{\vtop{
\eightpoint
\abovedisplayskip=6pt plus 3pt minus 3pt
\belowdisplayskip=6pt plus 3pt minus 3pt
\moreabstract\parindent=0 true in%
\caps{Abstract}: \ \ #1}}
\abovedisplayskip=12pt plus 3pt minus 9pt
\belowdisplayskip=12pt plus 3pt minus 9pt
\vskip 0.4 true in}
\def\moreabstract{%
\par \hsize = 5 true in \hangindent=0 true in \parindent=0.5 true in}

\def\caps#1{\smallcaps #1}
\def\nottencaps#1{\uppercase{#1}}

\def\firstbeginsection#1\par{\bigskip\vskip\parskip
\message{#1}\centerline{\caps{#1}}\nobreak\smallskip\noindent}

\def\beginsection#1\par{\vskip0pt plus.3\vsize\penalty-250
\vskip0pt plus-.3\vsize\bigskip\vskip\parskip
\message{#1}\centerline{\caps{#1}}\nobreak\smallskip\noindent}

\def\proclaim#1. #2\par{
\medbreak
\noindent{\caps{#1}.\enspace}{\it\rmpunctuation#2\par}
\ifdim\lastskip<\medskipamount \removelastskip
\penalty55\medskip\fi}

\def\Definition: #1\par{
\Head Definition: #1\par
\ifdim\lastskip<\medskipamount \removelastskip
\penalty55\medskip\fi}

\def\Problem #1: #2\par{
\Head Problem #1: #2\par
\ifdim\lastskip<\medskipamount \removelastskip
\penalty55\medskip\fi}

\def\sqr{\vcenter {\hrule height.3mm
\hbox {\vrule width.3mm height 2mm \kern2mm
\vrule width.3mm } \hrule height.3mm }}

\def\references#1{{
\frenchspacing
\eightpoint
\rmpunctuation
\halign{\bf##\hfil & \quad\vtop{\hsize=5.5 true
in\parindent=0pt\hangindent=3mm \strut\rm##\strut\smallskip}\cr#1}}}

\def\ref[#1]{{\bf [#1]}}

\catcode`@=11 
\def\vfootnote#1{\insert\footins\bgroup
\eightpoint
\interlinepenalty=\interfootnotelinepenalty
\splittopskip=\ht\strutbox
\splitmaxdepth=\dp\strutbox \floatingpenalty=20000
\leftskip=0pt \rightskip=0pt \spaceskip=0pt \xspaceskip=0pt
\textindent{#1}\footstrut\futurelet\next\fo@t}

\def\footremark{\insert\footins\bgroup
\eightpoint\it\rmpunctuation
\interlinepenalty=\interfootnotelinepenalty
\splittopskip=\ht\strutbox
\splitmaxdepth=\dp\strutbox \floatingpenalty=20000
\leftskip=0pt \rightskip=0pt \spaceskip=0pt \xspaceskip=0pt
\noindent\footstrut\futurelet\next\fo@t}
\catcode`@=12

\def\verses #1\par
{{\eightpoint\it\centerline{\hbox{\vbox{\halign{##\hfil\cr#1}}}}}}

\def\Bbb{\bf}
\def\E{{\Bbb E}}
\def\R{{\Bbb R}}
\def\Z{{\Bbb Z}}
\def\N{{\Bbb N}}
\def\C{{\Bbb C}}
\font\specialeightrm=cmr10 at 8pt
\def\R{\hbox{\rm I\kern-2pt R}}
\def\Z{\hbox{\rm Z\kern-3pt Z}}
\def\N{\hbox{\rm I\kern-2pt I\kern-3.1pt N}}
\def\C{\hbox{\rm \kern0.7pt\raise0.8pt\hbox{\specialeightrm I}\kern-4.2pt C}} 
\def\E{\hbox{\rm I\kern-2pt E}}

\def\Id{\hbox{\rm Id}}

\def\invc{{c^{-1}}}

\def\list#1,#2{#1_1$, $#1_2,\ldots,$\ $#1_{#2}}

\def\lnorm{\left\|}
\def\rnorm{\right\|}
\def\normo#1{\lnorm #1 \rnorm}

\def\lmod{\left|}
\def\rmod{\right|}
\def\modo#1{\lmod #1 \rmod}

\def\sign{\mathop{\hbox{\rm sign}}}


\def\tr{\mathop{\hbox{\rm tr}}}
\def\d{{\mathchoice{}{}{\hbox{\sevenrm d}}{\hbox{\fiverm d}}}}
\def\ad{{\mathchoice{}{}{\hbox{\sevenrm ad}}{\hbox{\fiverm ad}}}}


\noblankout

\setheadline The Distribution of Non-Commutative Rademacher Series\\
             Montgomery-Smith

{
\seventeenpoint
\centerline{The Distribution of Non-Commutative }
\smallskip
\centerline{Rademacher Series}
}
\bigskip\bigskip\medskip
\centerline{\caps{S.J.~Montgomery-Smith}%
\footnote{*}%
{Research supported in part by N.S.F.\ Grant DMS 9001796.}%
}
\smallskip
{
\eightpoint
\centerline{\it Department of Mathematics, University of Missouri,}
\centerline{\it Columbia, MO 65211.}
}
\bigskip\bigskip

\comment{
\verses
Cursed is the ground because of you;\cr
through painful toil you will eat of it\cr
all the days of your life.\cr
It will produce thorns and thistles for you,\cr
and you will eat the plants of the field.\cr
By the sweat of your brow\cr
you will eat your food\cr
until you return to the ground,\cr
since from it you were taken;\cr
for dust you are\cr
and to dust you will return.\cr
\hfill Genesis 3:17b--19 (N.I.V.)\cr
}

\verses
Cursed is the ground because of you;\cr
through painful toil you will eat of it\cr
all the days of your life.\cr
\hfill Genesis 3:17b (N.I.V.)\cr

\bigskip\bigskip

\Abstract

We give a formula for the tail of the distribution of the 
non-commutative Rademacher series, which generalizes the
result that is already available in the commutative case.  
As a result, we are able to calculate the norm
of these series in many rearrangement invariant spaces, generalizing
work of Pisier and Rodin and Semyonov.

\footremark{A.M.S.\ (1991) subject classification: 60-02,
60G50, 43A46, 46E30.}

\firstbeginsection 1.\ \ Introduction

The Rademacher functions are a sequence of independent 
random variables $r_n$\ such that 
$\Pr(r_n = \pm 1) = {1\over 2}$.  These functions have played
a very important role in mathematics, finding applications
in many parts of analysis, as well as other subjects like electronic
engineering.

One of the key inequalities concerning the Rademacher functions is due to
Khintchine
in 1923 \ref[Kh]:\ \ if $a_n$\ is a sequence of scalars, then for $0<p<\infty$\
$$ c_p \left(\sum_{n=1}^\infty \modo{a_n}^2\right)^{1/2}
   \le \left( \E \modo{ \sum_{n=1}^\infty a_n r_n }^p \right) ^{1/p}
   \le C_p \left(\sum_{n=1}^\infty \modo{a_n}^2\right)^{1/2} ,$$
where $C_p$\ and $c_p$\ are constants that depend upon $p$\ only. 
In particular, $C_p \le c\, \sqrt p$\ for $p\ge 1$.  (Throughout this
paper we will not be rigorous with infinite random sums --- an expression such
as the one above means that the random variable in the middle converges
in $L_p$\ if the right hand side is finite.)

Clearly, one would desire to find generalizations of such an important
inequality.  For example, one might like to calculate the norm of the
Rademacher series 
$\sum_{n=1}^\infty a_n r_n$\ in Orlicz or Lorentz spaces.
An obvious result (at least for real scalars) is the following:
$$ \normo{ \sum_{n=1}^\infty a_n r_n }_\infty
   = \sum_{n=1}^\infty \modo{a_n} .$$
However, another such generalization follows immediately from Khintchine's
inequality.  For a random variable $f$\
and $0<p<\infty$, let us denote by $\normo f_{\exp(t^p)}$\ the Orlicz
norm calculated using the Orlicz function  
$e^{t^p}-1$, i.e.
$$ \normo f_{\exp(t^p)} = \inf\bigl\{ \, \lambda :
   \E\bigl(\exp(\modo{f/\lambda}^p)\bigr) \le 2 \, \bigr\} .$$
Then
$$ \normo{ \sum_{n=1}^\infty a_n r_n }_{\exp(t^p)}
   \approx \left(\sum_{n=1}^\infty \modo{a_n}^2\right)^{1/2} ,$$
whenever $p\le 2$.  (Here, as in the rest of the paper, the expression
$A \approx B$\ means that $\invc A \le B \le c A$\ for some constant $c$.)

In 1975, Rodin and Semyonov \ref[R--S]  considered the value of the Rademacher 
series 
$\sum_{n=1}^\infty a_n r_n$\ in other rearrangement invariant spaces.
In particular, they showed that
$$ \normo{ \sum_{n=1}^\infty a_n r_n }_{\exp(t^p)}
   \approx \normo{(a_n)}_{q,\infty} ,$$
(see also \ref[P1]),
and that
$$ \normo{ \sum_{n=1}^\infty a_n r_n }_{\exp(t^p),r}
   \approx \normo{(a_n)}_{q,r} ,$$
whenever $p>2$, ${1\over p} + {1\over q} = 1$\ and $0<r<\infty$.
Here
$$ \eqalignno{
   \normo f_{\exp(t^p),r}
   &= \left( \int_0^1 \bigl(\log(1/t)\bigr)^{(r/p)-1} 
     \bigl(f^*(t)\bigr)^r \, {dt \over t}
     \right)^{1/r} , \cr
   \normo{(a_n)}_{q,r}
   &=
   \cases{ \left(\sum_{n=1}^\infty n^{(r/q)-1} a^*_n \right)^{1/r}
           & if $0<r<\infty$\cr
           \sup_{n\ge 1} n^{1/p} a^*_n
           & if $r=\infty$,\cr } \cr } $$
where $f^*$\ and $a_n^*$\ are the non-increasing rearrangements of $\modo{f}$\
and $\modo{a_n}$\ respectively.

In fact
they were able to show that if $X$\ is any symmetric sequence space
with Boyd indices strictly between $1$\ and $2$, then there exists a
rearrangement invariant space $Y$\ on probability space such that
$$ \normo{ \sum_{n=1}^\infty a_n r_n }_Y
   \approx \normo{(a_n)}_X .$$

There
still remained the question of finding tail distributions of Rademacher
series, that is, to find
$$ \Pr\left(\modo{\sum_{n=1}^\infty a_n r_n} > t \right)  $$
for every $t>0$.
This was answered in \ref[Mo] as follows.  Given a sequence $a = (a_n)$, we
will define its $K$-functional with respect to $\ell_1$\ and $\ell_2$\ to
be
$$ K_{1,2}(a,t) = \inf\left\{\,
   \normo{a'}_1 + t \normo{a''}_2 : a' + a'' = a \, \right\} .$$
These quantities play an important role in the theory of interpolation
of spaces (see \ref[B--S] or \ref[B--L]).  They are not so hard to
calculate, since there is the following formula due to Holmstedt \ref[Ho]:
$$ K_{1,2}(a,t) \approx
   \sum_{n=1}^{[t^2]} a_n^* + t \left(\sum_{n=[t^2]+1}^\infty (a_n^*)^2
   \right)^{1/2} .$$
Then we have the following results.
$$ \leqalignno{
   \Pr\left( \modo{\sum_{n=1}^\infty a_n r_n} > 
c \, K_{1,2}(t,a) \right) & \le c \, e^{-\invc t^2} , 
\cr
   \Pr\left( \modo{\sum_{n=1}^\infty a_n r_n} > 
\invc \, K_{1,2}(t,a) \right) & \ge \invc e^{-c t^2} . 
\cr } $$
(Here, as in the rest of the paper, the letter $c$\ denotes a positive
constant that changes with each occurrence.)
We remark that the hard part of this result, the lower bound, can
be deduced from a more general result contained in the book by 
Ledoux and Talagrand \ref[L--T], namely Theorem~4.15.  They obtain a 
Rademacher version of Sudakov's Theorem.

From this formula, and using known facts about the Hardy operators, it
is possible to reproduce all of the results of Rodin and Semyonov.
It is interesting to note that in order to obtain the lower bounds
of Rodin and Semyonov, one only 
requires the following estimate to be true:
$$  \Pr\left( \modo{\sum_{n=1}^\infty a_n r_n} > 
\invc \sum_{n=1}^{[t^2]} a_n^* \right)  \ge \invc e^{-c t^2} .$$
This bound has an extremely simple proof:\ \ simply consider the
random event $\{ r_{n_k} = \sign(a_{n_k}) \}$\ for an appropriate
sequence $\list n,{[t^2]}$.

We might also add that a consequence of the above result is the following.
If $t \le c\normo a_2 / \normo a_\infty$, then
$$ \Pr\left( \modo{\sum_{n=1}^\infty a_n r_n} > \invc t \normo a_2
   \right) \ge  \invc e^{-c t^2} .$$
This result can also be deduced from a result of Kolmogorov \ref[Ko]
(see also \ref[L--T] Chapter~4).

Recently, Hitczenko \ref[Hi] used the distribution
formula to obtain an asymptotically
more accurate version of Khintchine's original inequalities.  He showed that
$$ \left( \E \modo{ \sum_{n=1}^\infty a_n r_n }^p \right) ^{1/p}
   \approx K_{1,2}(a,\sqrt p) ,$$
for $p\ge 1$, where the constants of approximation do not depend upon $p$.

It has also been discovered that many of these results have vector valued
analogues (see \ref[D--M]).

\beginsection 2.\ \ The Non-Commutative Rademacher Series

Now we get to the main subject of this paper.  Non-commutative Rademacher
series arise in a natural way when one considers Fourier series on
non-commutative compact groups.  For example, a Sidon series on a 
non-commutative compact group has a distribution equivalent to a 
non-commutative Rademacher series (see \ref[F--R], \ref[H--R] and
\ref[A--M]).  They are also the natural things to consider if one
wishes to work with random Fourier series on non-commutative compact
groups (see \ref[M--P]).

They were considered by Fig\`a-Talamanca and 
Rider in \ref[F--R],
where they showed the non-commutative
analogue of the Khintchine inequalities (see also \ref[H--R]).  Many results
about them are also given in \ref[M--P].

We let $M_d$\ denote the vector space of $d$-dimensional matrices (i.e.\
$d\times d$\ matrices), and we let $O_d$\ denote the multiplicative
subgroup of orthogonal matrices.  Let $d_n$\ be a sequence of positive
integers, let $A_n$\ be a $d_n$-dimensional matrix, and let
$\epsilon_n$\ be a sequence of independent random variables such
that $\epsilon_n$\ takes values in $O_{d_n}$\ uniformly distributed
with respect to the Haar measure.  If $A$\ is a $d$-dimensional matrix,
we denote by $\tr(A)$\ the trace of $A$, that is, the sum of the diagonal
entries of $A$.

Then a non-commutative Rademacher series is a random variable of the following
form:
$$ S_\epsilon = \sum_{n=1}^\infty d_n \, \tr(\epsilon_n A_n) .$$

If $A$\ is a $d$-dimensional matrix, we define the singular values to be
the eigenvalues of $\sqrt{A^*A}$, where $A^*$\ is the transpose
of $A$.
We define the Schatten norms on $M_d$\ as follows:\enspace
if $A \in M_d$, set $\normo A_p$\ equal to the usual $\ell_p$\ sequence
norm of the singular values of $A$.  Thus $\normo A_\infty$\
is the usual operator norm of $A$\ on $d$-dimensional Hilbert space, 
$\normo A_2$\ is the Hilbert--Schmidt norm of $A$, and $\normo A_1$\ is
the trace class norm of $A$. 

The result of Fig\`a-Talamanca and Rider \ref[F--R] is the following:
$$ c_p \left(\sum_{n=1}^\infty d_n \normo{A_n}_2^2\right)^{1/2}
   \le
   \normo{ \sum_{n=1}^\infty d_n \, \tr(\epsilon_n A_n) }_p
   \le
   C_p \left(\sum_{n=1}^\infty d_n \normo{A_n}_2^2\right)^{1/2}, $$
for $0<p<\infty$.
Here $C_p \le c \sqrt p$\ for $p \ge 1$.  From this one can
obtain the result
$$ \normo{ \sum_{n=1}^\infty d_n \, \tr(\epsilon_n A_n) }_{\exp(t^2)}
   \approx
   \left(\sum_{n=1}^\infty d_n \normo{A_n}_2^2\right)^2 .$$
It is also true that
$$ \normo{ \sum_{n=1}^\infty d_n \, \tr(\epsilon_n A_n) }_\infty
   =
   \sum_{n=1}^\infty d_n \normo{A_n}_1 .$$

Let $s$\ denote the vector formed in the following manner.
First list the singular values of $A_n$, repeating each singular value
$d_n$\ times.
Combine these into one long list, rearranging them into decreasing order.
Then the above results can be written in the more suggestive form:
$$ \normo{S_\epsilon}_p \approx
   \normo{S_\epsilon}_{\exp(t^2)} \approx \normo s_2 
   \quad(0<p<\infty)
   \qquad\hbox{and}\qquad
   \normo{S_\epsilon}_\infty \approx \normo s_1 .$$

Pisier \ref[P1] was able to obtain partial non-commutative versions of the
results of Rodin and Semyonov.  He showed that
$$ \normo{S_\epsilon}_{\exp(t^p)} \le c \, \normo s_{q,\infty} ,$$
where $p>2$\ and ${1\over p} + {1\over q} = 1$.  He was not able to
obtain the lower bound.

The purpose of this paper is to show that all of these results for the
commutative Rademacher series also apply to the non-commutative case.  
The main
result is the following formulae for 
the distribution of the non-commutative Rademacher
series.

\proclaim Theorem 2.1.  The distribution of $S_\epsilon$\ is given by the 
following formulae.
$$ \leqalignno{
   \Pr( S_\epsilon > c \, K_{1,2}(t,s) ) & \le c \, e^{-\invc t^2} . 
   &\hbox{\rm \itemi{}}\cr
   \Pr( S_\epsilon > \invc \, K_{1,2}(t,s) ) & \ge \invc e^{-c t^2} . 
   &\hbox{\rm \itemii{}}\cr } $$

\proclaim Corollary 2.2.  We have the following for all $t > 0$:
$$ \invc \Pr \left( \sum_{n=1}^\infty s_n r_n > c t \right)
   \le
   \Pr \left( \sum_{n=1}^\infty d_n \, \tr(\epsilon_n A_n) > t \right) 
   \le
   c \, \Pr \left( \sum_{n=1}^\infty s_n r_n > \invc t \right) .$$

Now we are able to obtain the following results immediately from the
commutative case.

\proclaim Corollary 2.3.  We have the following inequalities.
\itemi $\normo{S_\epsilon}_{\exp(t^p)} \approx \normo s_{q,\infty}$\
for $p>2$\ and ${1\over p} + {1\over q} = 1$;
\itemii $\normo{S_\epsilon}_{\exp(t^p),r} \approx \normo s_{q,r}$\
for $p>2$, ${1\over p} + {1\over q} = 1$ and $0<r<\infty$;
\itemiii $\normo{S_\epsilon}_p \approx K_{1,2}(s,\sqrt p)$\
for $1\le p <\infty$\ with
with constants of approximation independent of $p$.
\moreproclaim

The proof of Theorem~2.1 is split into two halves.  In the first half
we make a number of changes to the problem, by showing that the problem
is equivalent to a similar result involving Gaussian matrices.

The second half contains the meat of the argument.
Part~(i) of Theorem~2.1 is essentially the result of Fig\`a-Talamanca and 
Rider combined
with a fairly straightforward interpolation argument.  It is part~(ii)
that provides the difficulties.  The argument proceeds by considering four
cases according to the nature of the sequence $s$.  

The arguments used in this paper assume that the matrices $A_n$\ are real
valued, but it is very easy to extend the results to the complex case as well.
Instead of using the non-commutative Rademacher functions, one should use
the non-commutative Steinhaus random variables, that is, $\xi_n$, where
$\xi_n$\ is uniformly distributed over the $d_n$-dimensional unitary matrices
with respect to Haar measure.  Then comparison results from \ref[M--P] combined
with Lemma~3.9 below will give the results.

This is probably a hard paper to read.  It certainly was a hard paper
to write.  As it says in Genesis 3:17, we eat of the ground
through painful toil.

\beginsection 3.\ \ The Proof of Theorem~2.1 --- Part I

The first 
observation is that
we may assume that all the matrices are diagonal
with entries from the non-negative reals.  This follows because
$A_n$\ may be factored $A_n = U_n D_n V_n$, where $U_n$\ and $V_n$\ are 
elements
of $O_{d_n}$, and $D_n$\ is diagonal with entries from the non-negative reals.  
But $\tr(\epsilon_n U_n D_n V_n) = \tr (V_n \epsilon_n U_n D_n)$, and
$V_n \epsilon_n U_n$\ has the same law as $\epsilon_n$.

Thus we will assume that 
$$ A_n = \left(\matrix{ a^n_1 & 0     & 0     & \cdots & 0         \cr
                        0     & a^n_2 & 0     & \cdots & 0         \cr
                        0     & 0     & a^n_3 & \cdots & 0         \cr
                        \vdots& \vdots& \vdots&        & \vdots    \cr
                        0     & 0     & 0     & \cdots & a^n_{d_n} \cr
         } \right) ,$$
where $\list a^n,{d_n}$\ are the singular values of $A_n$.

Let $G_n$\ be the matrix
$$ G_n = 
   {1\over \sqrt{d_n}}
   \left(\matrix{ g^n_{1,1} & g^n_{1,2} & g^n_{1,3} & \cdots & g^n_{1,d_n} \cr
                  g^n_{2,1} & g^n_{2,2} & g^n_{2,3} & \cdots & g^n_{2,d_n} \cr
                  g^n_{3,1} & g^n_{3,2} & g^n_{3,3} & \cdots & g^n_{3,d_n} \cr
                  \vdots    & \vdots    & \vdots    &        & \vdots      \cr
          g^n_{d_n,1} & g^n_{d_n,2} & g^n_{d_n,3} & \cdots & g^n_{d_n,d_n} \cr
         } \right) ,$$
where $(g^n_{i,j})$\ is a sequence of independent Gaussian random variables of
mean $0$\ and variance $1$.
We would like to compare $S_\epsilon$\ with the random variable
$$ S_G = \sum_{n=1}^\infty d_n \, \tr(G_n A_n) .$$
This random variable is particularly easy to understand --- it is simply
a Gaussian random variable:
$$ S_G = \sum_{n=1}^\infty \sqrt{d_n} \sum_{i=1}^{d_n}
         g^n_{i,i} a^n_i .$$
Unfortunately, this random variable is too large to give us the lower
bounds required for Theorem~2.1 part~(ii).  To get around this problem,
we split $G_n$\ as follows:
$$ G_n = G'_n + G''_n ,$$
where
$$ G'_n = G_n \chi_{\normo{G_n}_\infty \le \lambda}
   \hbox{\quad and\quad}
   G''_n = G_n \chi_{\normo{G_n}_\infty > \lambda} .$$
Here $\lambda$\ is a universal constant.  For the proof to work, $\lambda$\ 
needs 
to be sufficiently large.  As we proceed, we will make it clear where
the restrictions on $\lambda$\ are required.

We are also going to introduce the following random variables.  We let
$G^\d_n$\ denote the $d_n$-dimensional matrix consisting only of the diagonal
entries of $G_n$, and we let $G^\ad_n = G_n - G_n^\d$\ be the
$d_n$-dimensional matrix consisting of the off-diagonal
entries of $G_n$.  We can also split $G_n$\ in the following manner:
$$ G_n = G^{*}_n + G^{**}_n ,$$
where
$$ G^*_n = G_n (\chi_{\normo{G^\d_n}_\infty \le \lambda})
              (\chi_{\normo{G^\ad_n}_\infty \le \lambda})
   \hbox{\quad and\quad}
   G^{**}_n = G_n \chi_{\normo{G^\d_n}_\infty \vee \normo{G^\ad_n}_\infty
   > \lambda} .$$

The strategy will be to compare $S_\epsilon$\ with the random variables
$$ \eqalignno{
   S_{G'} &= \sum_{n=1}^\infty d_n \, \tr(G'_n A_n)
          = \sum_{n=1}^\infty d_n \, \tr(G_n A_n) 
            \chi_{\normo{G_n}_\infty \le \lambda} , \cr
  S_{G^*} &= \sum_{n=1}^\infty d_n \, \tr(G^*_n A_n)
          = \sum_{n=1}^\infty d_n \, \tr(G_n A_n) 
            \chi_{\normo{G^\d_n}_\infty \le \lambda}
            \chi_{\normo{G^\ad_n}_\infty \le \lambda} . \cr } $$

Now let us present the results that we will be requiring.  Note that we
denote the commutative Rademacher functions by $r_n$, so that they
will not be confused with the non-commutative Rademacher functions
$\epsilon_n$.  

The first pair of results we present are comparison
principles.  
Let us suppose that $V_n$\ is a sequence of random variables
taking values in $M_{d_n}$\ for which the sequence $(a_n V_n)$\ has the same
law as $(V_n)$\ for any $a_n \in O_{d_n}$.  We note that the random
variables $\epsilon_n$, $G_n$, $G'_n$\ and $G''_n$\
all have this property.  We also suppose that 
$x_n$\ is a sequence of $d_n$-dimensional matrices taking values
in a Banach space $B$.  
The first result is selected parts 
from \ref[M--P], Proposition~V.2.1.

\proclaim Lemma 3.1. Suppose that $T_n$\ is any element of $M_{d_n}$.
Then for $1\le p < \infty$\ we have that
$$ \left( \E \normo{\sum_{n=1}^\infty \tr(T_n V_n x_n)}^p
   \right)^{1/p}
   \le
   \sup_{n \ge 1} \normo{T_n}_\infty
   \left( \E \normo{\sum_{n=1}^\infty \tr(V_n x_n)}^p
   \right)^{1/p} .$$
We also have the following.
$$ \left( \E \normo{\sum_{n=1}^\infty \tr(\epsilon_n x_n)}^p
   \right)^{1/p}
   \le
   \sup_{n\ge1} \normo{\left(\E\modo{V_n}\right)^{-1}}_\infty
   \left( \E \normo{\sum_{n=1}^\infty \tr(V_n x_n)}^p
   \right)^{1/p} .$$

The next lemma is simply the commutative version of the same result,
and may be found in \ref[M--P], Theorem~4.9.

\proclaim Lemma 3.2.  Suppose that $v_n$\ is a sequence of independent,
real valued,
symmetric
random variables, and that $x_n$\ is a sequence of values from a Banach
space $B$.  Then for all $1 \le p < \infty$\ we have
$$ \left( \E \normo{\sum_{n=1}^\infty x_n r_n}^p \right)^{1/p}
   \le
   \sup_{n\ge 1} (\E\modo{v_n})^{-1}
   \left( \E \normo{\sum_{n=1}^\infty x_n v_n}^p \right)^{1/p} .$$

Next we present a couple of reflection principles.
The first will enable us to
remove some of the elements of $s$.  Let us suppose that
$V_n$\ is random variable taking its values in $M_{d_n}$\ such that
the sequence $(D_n V_n)$\ has the same law as $(V_n)$\ for any
sequence of diagonal matrices $D_n$\ whose diagonal entries are
$\pm 1$.  Notice that all the random variables we
have introduced have this property: $\epsilon_n$, $G_n$, $G'_n$,
$G''_n$, $G^*_n$\ and $G^{**}_n$.
Recall that we have supposed
that the matrices $A_n$\ are diagonal.

\proclaim Lemma~3.3.  Suppose that $A'_n$\ is a sequence of diagonal
matrices such that for each $n \ge 1$, each entry of $A'_n$\ is either
the same as the corresponding entry as $A_n$\ or it is $0$.  Then for
all $t>0$\ we have that
$$ \Pr \left(\sum_{n=1}^\infty d_n \, \tr(V_n A'_n) > t \right)
   \le 
   2\,\Pr \left(\sum_{n=1}^\infty d_n \, \tr(V_n A_n) > t \right) .$$

\Proof:  Notice that the random variables
$$ \sum_{n=1}^\infty d_n \, \tr(V_n A_n)
   \hbox{\quad and\quad}
   \sum_{n=1}^\infty d_n \, \tr(V_n(2A'_n-A_n)) $$
have the same law.  Thus
$$ \eqalignno{
   \Pr\left(\sum_{n=1}^\infty d_n \, \tr(V_n A'_n) > t \right)
   &\le \Pr\left(\sum_{n=1}^\infty d_n \, \tr(V_n A_n) > t \right)
   + \Pr\left(\sum_{n=1}^\infty d_n \, \tr(V_n (2A'_n - A_n)) > t \right) \cr
   & = 2\, \Pr \left(\sum_{n=1}^\infty d_n \, \tr(V_n A_n) > t \right) , \cr}$$
as required.
\endproof

The next lemma is simply the commutative version of the
above result, and is essentially the same
as \ref[Ka] Chapter~2 Theorem~5.

\proclaim Lemma 3.4.  Let $x_n$\ be any sequence of elements from
a Banach space $B$, and let $\alpha_n$\ be a sequence of values
taking only the values $0$\ or $1$.
Then for all $t>0$\ we have
$$ \Pr \left( \normo{\sum_{n=1}^\infty \alpha_n x_n r_n} > t \right)
   \le
   2 \, \Pr\left( \normo{\sum_{n=1}^\infty x_n r_n} > t \right) .$$

Now we present results concerning the behavior of the 
non-commutative
Gaussian random variables.

\proclaim Lemma 3.5.  For sufficiently large $\lambda$,
the following is true.
$$ \E \normo{G_n}_\infty \approx \E \normo{G'_n}_\infty \approx 1 . $$

\Proof: The statement $\E \normo{G_n}_\infty \approx 1$\ is given in 
\ref[M--P], Proposition~1.5.  That $\E \normo{G'_n}_\infty \approx 1$\
for sufficiently large $\lambda$\ then
follows by the monotone
convergence theorem.
\endproof

Let $\Id_n$\ denote the $d_n$-dimensional identity matrix.

\proclaim Lemma 3.6.  There exists constants $c_n$\
and $c'_n$\ that are uniformly
bounded above and below such that for sufficiently large $\lambda$\ we
have that
$$ \E \modo{G_n} = c_n \Id_n
   \hbox{\quad and\quad}
   \E \modo{G'_n} = c'_n \Id_n .$$

\Proof:  The first statement is from \ref[M--P], Corollary~1.8.  The second
statement has entirely the same proof.
\endproof

The next lemma uses a result
of C.~Borell \ref[Bo] (see also \ref[P2] or \ref[L--T]).

\proclaim Theorem~3.7.  Let $X$\ be 
a mean $0$\ Gaussian random variable taking values in a Banach space $B$.
Let 
$$ \sigma = \sup_{\normo{\phi}_{B^*}\le 1} \left(\E\modo{\phi(X)}^2 
   \right)^{1/2}.$$
Then for all $t>0$\ we have
$$ \Pr\left(
   \modo{\normo{X} - \E\normo{X}}
   \ge t \sigma \right)
   \le c\,e^{-\invc t^2} .$$

\proclaim Lemma~3.8. 
For $t$\ larger than some universal constant, we have that
$$ \Pr\left(\normo{G_n}_\infty > t \right) \le c\,e^{-\invc d_n t^2} .$$

\Proof:
Since $\E\normo{G_n}_\infty \approx 1$,
by Theorem~3.7, it is sufficient to show that if $\normo A_1 \le 1$, then
$$ \left(\E\modo{\tr(A^t G_n)}^2\right)^{1/2}
   \le {1\over\sqrt{d_n}} .$$
But 
$$ \tr( A^t G_n) = {1\over \sqrt{d_n}} \sum_{i=1}^{d_n} \sum_{j=1}^{d_n}
   g^n_{i,j} a^n_{i,j} ,$$
which is a Gaussian variable of variance
$$ {1\over d_n} \sum_{i=1}^{d_n} \sum_{j=1}^{d_n} (a^n_{i,j})^2
   = {1\over d_n} \normo A_2^2
   \le {1\over d_n} \normo A_1^2 .$$
\endproof

Now we present a principle from \ref[dP--M]
(see also \ref[A--M]) that allows us to obtain results about
distributions from $L_p$ norm results.

\proclaim Lemma 3.9.  Let $X$\ and $Y$\ be two random variables taking values
in the positive reals such that the following holds.
Whenever $X_m$\ and $Y_m$\ are independent random variables with the
same law as $X$\ and $Y$\ respectively, for all $M \in \N$\ we have that
$$ \leqalignno{
   \E \sup_{1\le m \le M} X_m 
   & \le c \, \E \sup_{1\le m \le M} Y_m ,
   &\hbox{\rm \itemi{}} \cr
   \left( \E \sup_{1\le m \le M} Y_m^2 \right)^{1/2} & 
   \le c \, \E \sup_{1\le m \le M} Y_m .
   &\hbox{\rm \itemii{}} \cr } $$
Then it follows that for all $t>0$ that
$$ \Pr (X > t) \le c \, \Pr (Y > \invc t) .$$

We can use this to prove a distributional comparison principle.

\proclaim Lemma 3.10.  Suppose that $v_n$\ is a sequence of real valued
symmetric independent
random variables, such that for any sequence of vectors $x_n$\ from a Banach
space $B$\ 
$$ \left( \E \normo{\sum_{n=1}^\infty x_n v_n}^2 \right)^{1/2}
   \le c \,
   \E \normo{\sum_{n=1}^\infty x_n v_n} .$$
Suppose also that
$$  \E\modo{v_n} \ge c^{-1} .$$
Then for any sequence of scalars $a_n$\ and for all $t > 0$\ we have
$$ \Pr\left( \modo{\sum_{n=1}^\infty a_n r_n} \ge t \right)
   \le c\,
   \Pr\left(\modo{\sum_{n=1}^\infty a_n v_n} \ge c^{-1} t \right) .$$

\Proof:  Let us set $X = \modo{\sum_{n=1}^\infty a_n r_n}$\ and
$Y = \modo{\sum_{n=1}^\infty a_n v_n}$.  
Let $r_{n,m}$\ be independent copies
of $r_n$\ and $v_{n,m}$\ be independent copies
of $v_n$\ for $1 \le m \le M$, 
and let $x_{n,m} \in \ell_\infty^M$\ be defined by
$$ x_{n,m} = 
   (0,0,\ldots,a_n,\ldots,0)
   \qquad\hbox{(the $a_n$\ is in the $m$th position)} .$$
Then notice that
$$ \eqalignno{
   \sup_{1\le m \le M} X_m
   & = 
   \normo{\sum_{m=1}^M \sum_{n=1}^\infty
   r_{n,m} x_{n,m}}_{\ell_\infty^M} ,\cr 
   \sup_{1\le m \le M} Y_m
   & = 
   \normo{\sum_{m=1}^M \sum_{n=1}^\infty d_n \, 
   v_{n,m} x_{n,m}}_{\ell_\infty^M} .\cr } $$
From Lemma~3.2, it then follows that for $p=1$, $2$\ that
$$ \left( \E \sup_{1\le m \le M} X_m^p \right)^{1/p}
   \le c\,
   \left( \E \sup_{1\le m \le M} Y_m^p \right)^{1/p} ,$$
and by hypothesis we have that
$$ \left( \E \sup_{1\le m \le M} Y_m^2 \right)^{1/2}
   \le c \,
   \E \sup_{1\le m \le M} Y_m .$$
Thus we may apply Lemma~3.9 and the result follows.
\endproof

The following lemma is an immediate corollary of \ref[M--P], Theorem~V.2.7.

\proclaim Lemma 3.11.  If $x_n$\ is a sequence of $d_n$-dimensional
matrices with entries in a Banach space $B$, then
$$ \E \normo{ \sum_{n=1}^\infty d_n \tr(\epsilon_n x_n) }
   \approx
   \left(\E \normo{ \sum_{n=1}^\infty d_n 
   \tr(\epsilon_n x_n) }^2 \right)^{1/2} . $$

Now we are ready to proceed with the main part of this section.  We are
going to use these results to show the following.

\proclaim Lemma 3.12.  For sufficiently large $\lambda$,
the following holds for all $t>0$:
$$ \eqalignno{
   \invc \Pr (S_{G'} > ct) &\le \Pr (S_\epsilon > t)
   \le c \, \Pr(S_{G'} > \invc t) , \cr
   \invc \Pr (S_{G^*} > ct) &\le \Pr (S_\epsilon > t)
   \le c \, \Pr(S_{G^*} > \invc t) . \cr } $$

First we will show the $L_p$-norm version of this result.

\proclaim Lemma 3.13. For sufficiently large $\lambda$, and
any $1\le p < \infty$\ we have that
$$ \left( \E \normo{\sum_{n=1}^\infty \tr(\epsilon_n x_n)}^p
   \right)^{1/p}
   \approx
   \left( \E \normo{\sum_{n=1}^\infty \tr(G'_n x_n)}^p
   \right)^{1/p} .$$

\Proof:  To show the left hand side is bounded by a constant times
the right hand side
is easy.  We apply the second part of Lemma~3.1 with $V_n = G'_n$,
using Lemma~3.6.

Next we show that the right hand side is bounded by a constant 
times the left hand
side.
Let us suppose that the random variables $G_n$\ are independent of the
random variables $\epsilon_n$.
Without loss of generality, we may suppose that 
measure space upon which the random variables exist is a product measure
of $\Omega_\epsilon$\ and $\Omega_G$, and that the random variables 
$\epsilon_n$\ depend only upon the $\Omega_\epsilon$\ co-ordinate, and
that the random variables $G_n$\ depend only upon the $\Omega_G$\ co-ordinate.
Let us denote integration with respect to the $\Omega_\epsilon$\ co-ordinate
by $\E_\epsilon$\ and 
integration with respect to the $\Omega_G$\ co-ordinate
by $\E_G$.

For each $\omega_G \in \Omega_G$, by Lemma~3.1, we have that
$$ \left( \E_\epsilon \normo{\sum_{n=1}^\infty 
   \tr(G'_n(\omega_G)\epsilon_n x_n)}^p
   \right)^{1/p}
   \le
   \sup_{n \ge 1} \normo{G'_n(\omega_G)}_\infty
   \left( \E_\epsilon \normo{\sum_{n=1}^\infty \tr(\epsilon_n x_n)}^p
   \right)^{1/p} .$$
Now we take $L_p$\ norms of both sides with respect with respect to the
$\Omega_G$\ co-ordinate.  We note that the sequence $(G'_n \epsilon_n)$\ has 
the same
joint law as $(G'_n)$, and that $\normo{G'_n}_\infty \le \lambda$.  Hence
we obtain that
$$ \left( \E \normo{\sum_{n=1}^\infty 
   \tr(G'_n  x_n)}^p
   \right)^{1/p}
   \le
   \lambda \,
   \left( \E \normo{\sum_{n=1}^\infty \tr(\epsilon_n x_n)}^p
   \right)^{1/p} ,$$
as desired.
\endproof

Now we need to be able to compare $S_{G'}$\ with $S_{G^*}$.  

\proclaim Lemma 3.14. Let $A$\ be a $d$-dimensional matrix, let $A^\d$\
be the matrix taking only the diagonal entries from $A$, and let $A^\ad$\
be the matrix taking only the non-diagonal entries from $A$, so that
$A = A^\d + A^\ad$.  Then
$$ \ds{1\over 2}
   \max\{\normo{A^\d}_\infty , \normo{A^\ad}_\infty \}
   \le \normo{A}_\infty
   \le 2 \, \max\{\normo{A^\d}_\infty , \normo{A^\ad}_\infty \} .$$

\Proof:  The right hand inequality follows immediately from the triangle
inequality.  To show the left hand inequality, note that
$$ \normo{A^\d}_\infty = \sup_{1\le i \le d} \modo{a_{i,i}}
   \le \normo A_\infty .$$
Finally, 
$$ \normo {A^\ad}_\infty \le \normo A_\infty + \normo {A^\d}_\infty
   \le 2 \normo A_\infty ,$$
as required.
\endproof

\proclaim Lemma 3.15. For all $t>0$\ we have
$$ \eqalignno{
   {1\over 2} \,
   \Pr\left( \normo{\sum_{n=1}^\infty \tr(G_n x_n) 
   \chi_{\normo{G_n}_\infty \le \lambda/2}r_n } > t \right)
   &\le
   \Pr\left( \normo{\sum_{n=1}^\infty \tr(\epsilon_n x_n)} > t \right) \cr
   &\le
   2 \,
   \Pr\left( \normo{\sum_{n=1}^\infty \tr(G_n x_n) 
   \chi_{\normo{G_n}_\infty \le 2\lambda}r_n } > t \right) . \cr } $$

\Proof:  As with the proof of Lemma~3.13,
we suppose that the random variables $G_n$\ are independent of the
random variables $r_n$.  We suppose that 
measure space upon which the random variables exist is a product measure
of $\Omega_r$\ and $\Omega_G$, and that the random variables 
$r_n$\ depend only upon the $\Omega_r$\ co-ordinate, and
that the random variables $G_n$\ depend only upon the $\Omega_G$\ co-ordinate.
Let us denote measure with respect to the $\Omega_r$\ co-ordinate
by $\Pr_r$.

By Lemma~3.14, the numbers
$$ { \chi_{\normo{G^\d_n}_\infty \le \lambda} 
     \chi_{\normo{G^\ad_n}_\infty \le \lambda}
     \over
     \chi_{\normo{G_n}_\infty \le 2 \lambda} }
\quad\hbox{and}\quad
   { \chi_{\normo{G_n}_\infty \le \lambda/2} 
     \over
     \chi_{\normo{G^\d_n}_\infty \le \lambda}
     \chi_{\normo{G^\ad_n}_\infty \le \lambda} } $$
take the values $0$\ or $1$.  Thus, by Lemma~3.4, it follows that
for each $\omega_G \in \Omega_G$\ and all $t>0$\ that
$$ \eqalignno{
   & {1\over 2} \,
   \Pr\nolimits_r \left( \normo{\sum_{n=1}^\infty \tr(G_n(\omega_G) x_n) 
   (\chi_{\normo{G_n(\omega_G)}_\infty \le \lambda/2})r_n }
   >t \right) \cr
   &\le c\,
   \Pr\nolimits_r \left( \normo{\sum_{n=1}^\infty \tr(G^*_n(\omega_G) x_n) r_n }
   >t \right) \cr
   &\le 2\,
   \Pr\nolimits_r \left( \normo{\sum_{n=1}^\infty \tr(G_n(\omega_G) x_n) 
   (\chi_{\normo{G_n(\omega_G)}_\infty \le 2\lambda})r_n }
   > t \right) .\cr}$$
Now, taking expectations on both sides with respect to $\Omega_G$, the
result follows.
\endproof

Now we are ready to combine these results.

\Proofof Lemma 3.12:  We first prove the first inequality.
In order to apply Lemma~3.9, let us set $X = \modo{S_{G'}}$\ and
$Y = \modo{S_\epsilon}$.  Let $\epsilon_{n,m}$\ be independent copies
of $\epsilon_n$\ and $G_{n,m}$\ be independent copies
of $G_n$\ for $1 \le m \le M$, 
and let $x_{n,m}$\ be 
diagonal matrices with diagonal entries in $\ell_\infty^M$:
$$ x^{n,m}_{i} = 
   (0,0,\ldots,a^n_{i},\ldots,0)
   \qquad\hbox{(the $a^n_{i}$\ is in the $m$th position)} .$$
Then notice that
$$ \eqalignno{
   \sup_{1\le m \le M} X_m
   & = 
   \normo{\sum_{m=1}^M \sum_{n=1}^\infty d_n \, 
   \tr(G'_{n,m} x_{n,m})}_{\ell_\infty^M} ,\cr 
   \sup_{1\le m \le M} Y_m
   & = 
   \normo{\sum_{m=1}^M \sum_{n=1}^\infty d_n \, 
   \tr(\epsilon_{n,m} x_{n,m})}_{\ell_\infty^M} .\cr } $$
From Lemma~3.13, it then follows that for $p=1$, $2$\ that
$$ \left( \E \sup_{1\le m \le M} X_m^p \right)^{1/p}
   \approx
   \left( \E \sup_{1\le m \le M} Y_m^p \right)^{1/p} ,$$
and from Lemma~3.11, we have that
$$ \left( \E \sup_{1\le m \le M} Y_m^2 \right)^{1/2}
   \approx
   \E \sup_{1\le m \le M} Y_m .$$
Thus, we also have that
$$ \left( \E \sup_{1\le m \le M} X_m^2 \right)^{1/2}
   \approx
   \E \sup_{1\le m \le M} X_m .$$
Thus we may apply Lemma~3.9 twice, once with the roles of $X_m$\ and
$Y_m$\ reversed, and the result follows.

The second inequality now follows from Lemma~3.15.
\endproof

\beginsection 4.\ \ The Proof of Theorem~2.1 --- Part II

We will first show the first half of Theorem~2.1.

\proclaim Proposition 4.1. The following is true for all $t>0$.
$$ \Pr ( S_\epsilon > c \, K_{1,2}(t,s) ) \le c \, e^{-\invc t^2} .$$

\Proof:  Choose sequences $s'$ and $s''$\ such that $s = s' + s''$\ and
$$ K_{1,2}(t,s) \ge \ts{1\over 2}(\normo{s'}_1 + t \normo{s''}_2) .$$
We may assume that if a certain number occurs several times in the sequence
$s$, then for each occurrence this number is split identically between
$s'$\ and $s''$.  Thus we may know that there exist
sequences of matrices $(A_n')$\ and
$(A_n'')$\ such that $A_n = A_n' + A_n''$\ and such $s'$\ comes from repeating
$d_n$\ times
the singular values of $A'_n$, and $s''$\ comes from repeating
$d_n$\ times
the singular values of $A''_n$.  

From the result of Fig\`a-Talamanca and 
Rider, we know that
$$ \Pr \left( \sum_{n=1}^\infty d_n \tr(\epsilon_n A''_n) 
              > ct \normo{s''}_2 \right) 
   \le {\sqrt{p^p} \over t^p} \le c \, e^{-\invc t^2} .$$
(Here we chose $p = t^2/2$).
It is also clearly evident that 
$$ \sum_{n=1}^\infty d_n \tr(\epsilon_n A'_n) 
   \le \sum_{n=1}^\infty d_n \, \normo{A'_n}_1 = \normo{s'}_1 .$$
Thus
$$ \eqalignno{
   \Pr ( S_\epsilon > 2c \, K_{1,2}(t,s) )
   &\le \Pr \left( \sum_{n=1}^\infty d_n\tr(\epsilon_n A'_n)
              + \sum_{n=1}^\infty d_n\tr(\epsilon_n A''_n) 
              > c(\normo{s'}_1 + t \normo{s''}_2) \right) \cr
   &\le \Pr \left( \sum_{n=1}^\infty d_n\tr(\epsilon_n A''_n) 
              > ct \normo{s''}_2 \right) \cr
   &\le c \, e^{-\invc t^2} ,\cr }$$
and the result follows.
\endproof

Now we finally come to the hard part of this paper:\enspace to show the 
second part of Theorem~2.1.  We will proceed by considering
three cases.  All of the arguments will make heavy use of the
approximation 
$$ \invc e^{-ct^2} \le \Pr(g>t) \le c e^{-\invc t^2} 
   \qquad (t > 0),$$
whenever $g$\ is a Gaussian random variable of mean $0$\ and variance $1$.

Our use of the letter $c$\ becomes confusing at this point.  Thus from now
on we will use subscripts on the letter $c$\ to denote different values.
However, the same subscripted letter $c$\ may take different values from
result to result and proof to proof.

The first case will be dealt with by the following 
result.

\proclaim Proposition 4.2.  For sufficiently large $\lambda$,
the exist numbers $c_1$\ and $c_2$\ such that for all
integers $t \ge 1$.
$$ \Pr ( S_{G^*} > c_1^{-1} \sum_{m=1}^t s_m) \ge e^{-c_2 t} .$$

\Proof:  We suppose that $\list s,t$\ is made up
as follows:\enspace for each $n\ge 1$\ and $1\le i \le d_n$, we pick
$0 \le K_{n,i} \le d_n$.  Then the sequence 
$(\list s,t)$\ consists of the $a^n_i$, each one repeated $K_{n,i}$\ times.
Let $L$\ be the number of pairs $(n,i)$\ such that we have
$K_{n,i} \ne 0$.  Define the following events.
$$ \eqalignno{
   B_{n,i} &= \Bigl\{\,\sqrt{d_n} g^n_{i,i} a_i^n \ge c_1^{-1} K_{n,i} a^n_i
   \,\Bigr\} ,\cr
   C_{n} &= \Bigl\{\,\sqrt{d_n}\normo{G^\d_n}_\infty = 
   \sup_{1\le i\le d_n} \modo{g^n_{i,i}} \le \sqrt{d_n}\lambda\,\Bigr\} ,\cr
   D_{n} &= \Bigl\{\,\normo{G^\ad_n}_\infty \le \lambda\,\Bigr\} .\cr } $$
By Lemma~3.3, we are really asking for a lower bound for
the probability of the event
$$ \sum_{n=1}^\infty \sum_{i=1}^{d_n}
   \sqrt{d_n} g^n_{i,i} a^n_i
   (\chi_{\normo{G^\d_n}_\infty \le \lambda})
   (\chi_{\normo{G^\ad_n}_\infty \le \lambda})
   \ge
   c_1^{-1} \sum_{n=1}^\infty \sum_{i=1}^{d_n} K_{n,i} a_i^n  .$$
However, we notice that this event contains 
$$ \bigcap_{(n,i):K_{n,i} \ne 0} B_{n,i} \cap C_{n}
   \cap D_{n} .$$
\noindent
Now, 
$$ \bigcap_{(n,i):K_{n,i} \ne 0}
   B_{n,i} \cap C_{n} = 
   \bigcap_{(n,i):K_{n,i} \ne 0}
   \Bigl\{\, c_1^{-1} K_{n,i} / \sqrt{d_n} \le
   g^n_{i,i} \le \sqrt{d_n} \lambda \,\Bigr\} .$$
Since
$K_{n,i} \le d_n$, if $\lambda$\ and $c_1$\ are chosen large enough, then
we see that
$$ \Pr \left( \bigcap_{(n,i):K_{n,i} \ne 0} B_{n,i} \cap C_{n} \right)
   \ge 
   c_3^{-L} \exp\left(-c_4 \sum_{n=1}^\infty \sum_{i=1}^{d_n} K_{n,i}^2 / d_n 
   \right) .$$
Event $D_{n}$ is independent of $B_{n,i}$ and $C_{n}$.  
By Lemma~3.8, it follows
that for each number $n$\ that
$$ \Pr( \normo{G^\ad_n}_\infty \le \lambda )
   \ge \Pr( \normo{G_n}_\infty \le \lambda /2 ) 
   \ge c_5^{-1}. $$
Thus 
$$ \Pr \left( \bigcap_{(n,i):K_{n,i} \ne 0} B_{n,i} \cap C_{n} \cap D_{n}
   \right) 
   \ge 
   (c_3 c_5)^{-L} \exp\left(-c_4\sum_{n=1}^\infty 
   \sum_{i=1}^{d_n} K_{n,i}^2 / d_n 
   \right) .$$
Now $K_{n,i}^2 / d_n \le K_{n,i}$, and further, if $u \ge 1$, then
$(c_3 c_5)^{-1} e^{-c_4 u} \ge e^{-c_2 u}$.  Hence the probability that 
we require
is bounded below by
$$ \eqalignno{
   (c_3 c_5)^{-L} \exp\left(-c_4\sum_{n=1}^\infty 
   \sum_{i=1}^{d_n} K_{n,i}
   \right)
   &= 
   \prod_{(n,i):K_{n,i}\ne0} \left( (c_3 c_5)^{-1} \exp(-c_4 K_{n,i}) \right)
   \cr
   &\ge
   \prod_{(n,i):K_{n,i}\ne0} \exp(-c_2 K_{n,i}) \cr
   &=
   \exp\left(-c_2\,\sum_{n=1}^\infty \sum_{i=1}^{d_n} K_{n,i} \right) \cr
   &= 
   e^{-c_2t} ,\cr}$$
as desired.
\endproof

Now we are ready for the second case.

\proclaim Proposition 4.3. Fix $t > 0$.
Suppose that there is a number $c_1$\
such that for all $n\ge 1$\ and $1\le i \le d_n$
that either $a_i^n = 0$\ or 
$$ c_1^{-1} \normo s_2 \le t \sqrt{d_n} a_i^n \le c_1\sqrt{d_n}\,\normo s_2 .$$
Then for sufficiently large $\lambda$, there are numbers $c_2$\ and
$c_3$, depending only on $c_1$\ and $\lambda$, such that
$$ \Pr( S_{G^*} \ge c_2^{-1} t \normo s_2) \ge e^{-c_3 t^2} .$$

\Proof:  The second case has a very similar proof to the first case.
First, without loss of generality, we may suppose that $A_n \ne 0$\ for
all $n \ge 1$.  Also recall that
$$ \normo s_2 = \left( \sum_{n=1}^\infty \sum_{i=1}^{d_n} d_n (a_i^n)^2
                \right)^{1/2} .$$
Define the events
$$ \eqalignno{
   B_{n,i} &= \cases{
     \Bigl\{\, c_2^{-1} t a^n_i \sqrt{d_n} / \normo s_2 
     \le g^n_{i,i} \le \sqrt{d_n} \lambda \,\Bigr\}
     & if $a_i^n \ne 0$,\cr
     &\cr
     \Bigl\{\, g^n_{i,i} \le \sqrt{d_n} \lambda \,\Bigr\}
     & if $a_i^n = 0$,\cr}\cr
   C_{n} &= \Bigl\{\, \normo{G^\ad_n}_\infty \le \lambda \,\Bigr\} . \cr } $$
By Lemma~3.3, we are looking for a lower bound for the event
$$ \sum_{n=1}^\infty \sum_{i=1}^{d_n}
   \sqrt{d_n} g^n_{i,i} a^n_i
   \chi_{\normo{G^\d_n}_\infty \le \lambda}
   \chi_{\normo{G^\ad_n}_\infty \le \lambda}
   \ge
   c_2^{-1} t \normo s_2 ,$$
This event contains 
$$ \bigcap_{(n,i)} B_{n,i} \cap C_n .$$
Let us first consider $\Pr(B_{n,i})$\ in the case when $a_i^n \ne 0$.
Since 
$t \sqrt{d_n} a_i^n \le c_1\sqrt{d_n}\,\normo s_2$, we
see that for sufficiently large $\lambda$\
that $2 c_2^{-1} t a^n_i \sqrt{d_n} / \normo s_2 \le \sqrt{d_n} \lambda$,
and hence
$$ \Pr(B_{n,i}) \ge c_4^{-1} 
   \exp\left( - c_5 t^2\, {d_n (a_i^n)^2 \over \normo s_2^2} \right) .$$
Now $c_1^{-1} \normo s_2 \le t \sqrt{d_n} a_i^n$, and if $u > c_1^{-2}$,
then $c_4 e^{-c_5 u} \ge e^{-c_6 u}$, and hence
$$ \Pr(B_{n,i}) \ge
   \exp\left( - c_6 t^2\, {d_n (a_i^n)^2 \over \normo s_2^2} \right) $$
Hence for each $n \ge 1$
$$ \Pr\left(\bigcap_{i:a_i^n \ne 0} B_{n,i}\right) \ge 
   \exp\left( - c_6 t^2 \sum_{i=1}^{d_n}
   {d_n (a_i^n)^2 \over \normo s_2^2} \right) .$$
Also
$$ \Pr\left(C_n \cap \bigcap_{i:a_i^n = 0} B_{n,i} \right) \ge
   \Pr(\normo{G_n} \le \lambda/2)
   \ge c_7^{-1},$$ 
where the last inequality follows from Lemma~3.8
if $\lambda$\ is sufficiently
large.  Since $c_7 e^{-c_6 u} \ge e^{-c_3 u}$\ whenever $u > c_1^{-2}$, 
it follows that 
$$ \Pr\left(C_n \cap \bigcap_i B_{n,i} \right) \ge
   \exp\left( - c_3 t^2 \sum_{i=1}^{d_n}
   {d_n (a_i^n)^2 \over \normo s_2^2} \right) .$$
Hence, 
$$ \Pr\left(\bigcap_{(n,i)} B_{n,i} \cap C_n \right) \ge 
   \exp\left( - c_3 t^2 \sum_{n=1}^\infty \sum_{i=1}^{d_n}
   {d_n (a_i^n)^2 \over \normo s_2^2} \right)  = e^{-c_3t^2} ,$$
as desired.
\endproof

Now for the third case.  The argument that follows
was suggested by the proof of Proposition~4.13 in \ref[L--T]. 

\proclaim Proposition 4.4. Fix $t > 0$.
Suppose that there is a number $c_1$\
such that for all $n\ge 1$\ and $1\le i \le d_n$
that either $a_i^n = 0$\ or 
$$ c_1^{-1} \normo s_2 /\sqrt{d_n} \le
   t \sqrt{d_n} a_i^n \le c_1 \normo s_2 .$$
Then for sufficiently large $\lambda$, 
there are numbers $c_2$\ and $c_3$, depending only on $c_1$\ and $\lambda$,
such that
$$ \Pr( S_{G'} \ge t \normo s_2) \ge c_2^{-1} e^{-c_3t^2} .$$

\Proof:  First note that for $u>0$\ that
$$ \eqalignno{
   \Pr\left( d_n \tr(A_n G_n) \chi_{\normo{G_n}_\infty > \lambda} > u \right)
   &\le \min\left\{ \Pr\left( d_n \tr(A_n G_n) > u \right) , \,
                   \Pr\left( \normo{G_n}_\infty > \lambda \right) \right\} \cr
   &\le c_4 \min\left\{ 
              \exp\left(-{c_4^{-1} u^2\over d_n \normo{A_n}_2^2}\right)
              ,\, \exp(-c_4^{-1} \lambda^2 d_n) \right\} . \cr } $$
(Here we used Lemma~3.8.)
Now let 
$$ \theta = {t\lambda^{1/2}\over \normo s_2} .$$
Since $ c_1^{-1} \normo s_2 /\sqrt{d_n} \le
t \sqrt{d_n} a_i^n \le c_1\normo s_2$\ whenever $a_i^n \ne 0$, it follows that
$$ c_1^{-2} {\normo s_2^2 \over d_n}
   \le t^2 d_n \sum_{i : a_i^n \ne 0} (a_i^n)^2
   \le c_1^2 d_n \normo s_2^2 ,$$
i.e., $ c_1^{-1} \normo s_2/d_n \le t\,\normo{A_n}_2
\le c_1 \normo s_2$. Hence
$$ c_1^{-1}\lambda^{1/2} / d_n \le \theta \normo{A_n}_2
   \le c_1\lambda^{1/2} .$$
Now
$$ \eqalignno{
   &\E \left( \exp\left(\theta d_n \tr(A_n G_n) 
   \chi_{\normo{G_n}_\infty > \lambda}
   \right) \right)
   \le 1 + \int_0^\infty \theta e^{\theta u} 
   \Pr \left( d_n \tr(A_n G_n) 
   \chi_{\normo{G_n}_\infty > \lambda} > u \right) \, du \cr
   &\le
   1 + c_4 \int_{0}^{\lambda d_n \normo{A_n}_2 }
   \theta e^{\theta u} \exp(-c_4^{-1} \lambda^2 d_n) \, du
   + c_4 \int_{\lambda d_n \normo{A_n}_2 }^\infty
   \theta e^{\theta u} \exp\left(-{c_4^{-1} u^2\over d_n \normo{A_n}_2^2}\right)
   \, du . \cr }$$
Furthermore,
$$ \eqalignno{
   \int_{\lambda d_n \normo{A_n}_2 }^\infty
     \theta e^{\theta u} 
     \exp\left(-{c_4^{-1} u^2\over d_n\normo{A_n}_2^2}\right)
     \, du
   &\le
   \int_{\lambda d_n \normo{A_n}_2 }^\infty
     \theta e^{\theta u} \exp\left(-{c_4^{-1} \lambda u \over \normo{A_n}_2
     }\right) \, du \cr 
   & = {\theta \normo{A_n}_2\over c_4^{-1}\lambda - \theta \normo{A_n}_2}
       \exp\left( \left(\theta \normo{A_n}_2- c_4^{-1}\lambda\right) \lambda
       d_n \right) \cr
   & \le c_5 \exp(-c_6^{-1} \lambda^2 d_n) ,\cr}$$
when $\lambda$\ is sufficiently large, because $\theta \normo{A_n}_2
\le c_1 \lambda^{1/2}$.
Similarly
$$ \eqalignno{
   \int_{0}^{\lambda d_n \normo{A_n}_2 }
     \theta e^{\theta u} \exp(-c_4^{-1} \lambda^2 d_n) \, du 
   & = \exp\left(-c_4^{-1} \lambda^2 d_n + \lambda \theta d_n \normo{A_n}_2
       \right) \cr
   &\le \exp(-c_6^{-1} \lambda^2 d_n) ,\cr}$$
when $\lambda$\ is sufficiently large.

Now, 
$$ 1 \le  c_1^2 \lambda^{-1} \theta^2 d_n^{2} \normo{A_n}_2^2 ,$$
and since $e^{-u} \le 1/u$\ for $u > 0$, 
$$ \exp(-c_6^{-1} \lambda^2 d_n) \le c_6 \lambda^{-2} d_n^{-1} ,$$
and so
$$ \exp(-c_6^{-1} \lambda^2 d_n) 
   \le c_1^2c_6 \lambda^{-3} \theta^2 d_n \normo{A_n}_2^2 .$$
Hence
$$ \eqalignno{
   \E \left( \exp\left(\theta d_n \tr(A_n G_n) 
   \chi_{\normo{G_n}_\infty > \lambda}
   \right) \right)
   &\le 1 + c_4(1+c_5)c_1^2c_6 \lambda^{-3} \theta^2 d_n \normo{A_n}_2^2 \cr
   & \le \exp(c_7 \lambda^{-3} \theta^2 d_n \normo{A_n}_2^2) .\cr } $$
Now, 
$$ \eqalignno{
   \E (\exp(\theta S_{G''}))
   &= \E \left(\prod_{n=1}^\infty \exp\bigl(\theta d_n \tr(A_n G''_n)\bigr)
      \right) \cr
   &= \prod_{n=1}^\infty \E \Bigl( \exp\bigl(\theta d_n \tr(A_n G''_n)\bigr)
      \Bigr) \cr 
   &\le \prod_{n=1}^\infty \exp(c_7 \lambda^{-3} \theta^2 d_n \normo{A_n}_2^2) 
        \cr
   &= \exp(c_7 \lambda^{-3} \theta^2 \normo s_2^2 ) .\cr}$$
So,
$$ \Pr\left( \exp(\theta S_{G''}) > 
      \exp(c_7 \lambda^{-3} \theta^2 \normo s_2^2 
      + \lambda^{1/2} t^2) \right) \le e^{-\lambda^{1/2} t^2} ,$$
that is,
$$ \Pr\left( S_{G''} > (1+c_7 \lambda^{-5/2}) t \normo s_2 \right) 
   \le e^{-\lambda^{1/2} t^2} .$$
If $\lambda > c_7^{-2/5}$, then
$$ \Pr\left( S_{G''} > 2 t \normo s_2 \right) 
   \le e^{-\lambda^{1/2} t^2} .$$
To finish, we note that
$$ \Pr\left( S_{G'} >  t \normo s_2 \right)
   \ge \Pr\left( S_{G} >  2t \normo s_2 \right)
       - \Pr\left( S_{G''} >  t \normo s_2 \right)
   \ge c_8^{-1} e^{-c_8t^2} - e^{-\lambda^{1/2} t^2/4} .$$
Thus if $\lambda$\ is sufficiently large, then
$$ \Pr\left( S_{G'} >  t \normo s_2 \right) \ge c_2^{-1} e^{-c_3t^2} $$
for $ t > {1\over 2}$.  

If $t \le {1\over 2}$, then we can use the following inequality 
(see \ref[Ka], Chapter~1):\ \ if
$X$\ is a positive random variable, then
$$ \Pr(X \ge \normo X_1/2) \ge {\normo X_1^2 \over 4 \normo X_2^2 } .$$
Take $X = \modo{S_\epsilon}^2$.  By the result of 
Fig\`a-Talamanca and Rider, and Lemma~3.13 it follows that $\normo X_2 \le
c_9 \normo X_1$, and the result follows.
\endproof

The fourth case follows by comparing the non-commutative case with
the commutative case.

\proclaim Proposition 4.5. Fix $t>0$.
Suppose that there is a number
$c_1$\ such that for all $n\ge 1$\ and $1\le i \le d_n$
$$ t \sqrt{d_n} a_i^n \le  c_1 \normo s_2 / \sqrt{d_n} .$$
Then there is a number $c_2$, depending only on $c_1$, such that
$$ \Pr( S_{\epsilon} \ge c_2^{-1} t \normo s_2) \ge c_2^{-1} e^{-c_2t^2} .$$

\Proof:  We apply Lemma~3.10 with 
$$ v_n = {\sqrt{d_n} \tr(A_n \epsilon_n) \over  \normo{A_n}_2} $$
and $a_n = \sqrt{d_n} \normo{A_n}_2$\ to deduce that
$$ \Pr\left( S_{\epsilon} >  t \normo s_2 \right)
   \ge c_3^{-1}
   \Pr\left( \sum_{n=1}^\infty
   \sqrt{d_n} \normo{A_n}_2 r_n >  c_3 t \normo s_2 \right) .$$
From the hypothesis, we have that 
$$ t \sqrt{d_n} \normo{A_n}_2 \le c_1 \normo s_2 ,$$
and hence by the commutative Rademacher series result, it follows that
$$ \Pr\left( \sum_{n=1}^\infty
   \sqrt{d_n} \normo{A_n}_2 r_n >  c_4^{-1} t \normo s_2 \right)
   \ge c_4^{-1} e^{-c_4t^2} ,$$
as required.
\endproof

Now we are finally ready to put the pieces together.  Let us restate
the theorem we are attempting to prove.

\proclaim Theorem 2.1.  The distribution of $S_\epsilon$\ is given by the 
following formulae.
$$ \leqalignno{
   \Pr( S_\epsilon > c \, K_{1,2}(t,s) ) & \le c \, e^{-\invc t^2} . 
   &\hbox{\rm \itemi{}}\cr
   \Pr( S_\epsilon > \invc \, K_{1,2}(t,s) ) & \ge \invc e^{-c t^2} . 
   &\hbox{\rm \itemii{}}\cr } $$

\Proof:  Part~(i) is simply Proposition~4.1.  To prove part~(ii),
we may suppose that $t \ge 1$.
Now, note that we have the following bound:
$$ K_{1,2}(t,s) \le
   \sum_{m=1}^{[t^2]} s_m
   + t \, \left( \sum_{m=[t^2]+1}^\infty (s_m)^2 \right)^{1/2} .$$
Then we have
two possibilities.

\Head Case 1: The first possibility is that
$$ \sum_{m=1}^{[t^2]} s_m \ge \ts {1\over 2} \, K_{1,2}(t,s) .$$
In that case, the result follows by Lemma~3.12 and Proposition~4.2.

\Head Case 2: Otherwise, we know that 
$$ t \, \left( \sum_{m=[t^2]+1}^\infty (s_m)^2 \right)^{1/2}
   \ge \ts{1\over 2} \, K_{1,2}(t,s) .$$
We also know that
$$ t \, \left( \sum_{m=[t^2]+1}^\infty (s_m)^2 \right)^{1/2}
   \ge \sum_{m=1}^{[t^2]} s_m \ge [t^2] s_{[t^2]} ,$$
since the sequence $(s_m)$\ is in decreasing order.
Hence, if $m \ge [t^2]$, we have that
$$ 2t s_m \le \left( \sum_{m=[t^2]+1}^\infty (s_m)^2 \right)^{1/2} .$$
Let $M$\ be the least number $m$\ such that $s_m = s_{[t^2]+1}$.
Let us replace the matrices $A_n$\ with matrices $A'_n$\ that drop the entries
that correspond to $s_m$\ for $m < M$.  Thus the new sequence $s'$\ formed
satisfies the following.
$$ \normo{s'}_2 \ge \ts{1\over 2} \, K_{1,2}(t,s)
   \hbox{\quad and\quad}
   2t s'_m \le \normo{s'}_2 .$$
Thus, we may replace the matrices $A_n$\ with $A'_n$, and,
using Lemma~3.3, we are reduced to
showing the following:\enspace subject to the restriction that
$$ 2t a_i^n \le \normo s_2 ,$$
we desire to show that
$$ \Pr ( S_\epsilon \ge 2 \invc t \normo s_2 ) \ge 2 \invc e^{-ct^2} .$$
To prove this, we will split the entries of the matrices into three parts.
Let 
$$ \eqalignno{
   B_1 &= \{\, (n,i) : 
   \normo s_2 < t \sqrt{d_n} a_i^n \, 
   \} , \cr
   B_2 &= \{\, (n,i) : 
   \normo s_2 / \sqrt{d_n} < t \sqrt{d_n} a_i^n \le \normo s_2 \, 
   \} , \cr
   B_3 &= \{\, (n,i) : 
   t \sqrt{d_n} a_i^n \le \normo s_2 / \sqrt{d_n} \, 
   \} . \cr } $$
Then for one of $j = 1$, $2$, $3$, we have that
$$ \sum_{(n,i) \in B_j} d_n (a_i^n)^2 \ge \ts{1\over3} \normo s_2^2 .$$
In that case, we can replace the matrices $A_n$\ with matrices that only
take those entries that are in the set $B_j$.  Now the result follows
by Lemma~3.12, Lemma~3.3 and Proposition~4.3, 4.4 or~4.5.

\endproof

\beginsection Acknowledgements

The author would like to thank Nakhl\'e Asmar for both useful 
discussions and warm friendship while this paper was being prepared.
He would also like to express gratitude to the referee for useful
corrections and comments.

\beginsection References

\references{
A--M & N.~Asmar and S.J.~Montgomery-Smith,\rm\ On the distribution of 
Banach valued
Sidon spectral functions,\sl\ Arkiv Mat.\ (to appear).\cr
B--S & C.~Bennett and R.~Sharpley,\sl\ Interpolation of Operators,\rm\
Academic Press.\cr
B--L & J.~Bergh and J. L\"ofstr\"om,\sl\ Interpolation Spaces,\rm\
Springer-Verlag, 1976.\cr
Bo & C.~Borell,\rm\ The Brunn--Minkowski inequality in Gauss Space,\sl\
Invent.\
Math.\ {\bf 30} (1975) 207--216.\cr
dP--M & V.H.~de~la~Pe\~na and S.J.~Montgomery-Smith,\rm\  
Decoupling inequalities for tail 
probabilities 
of multilinear forms of symmetric and hypercontractive variables,\sl\ 
submitted.\cr
D--M. & S.J.~Dilworth and S.J.~Montgomery-Smith,\rm\ 
The distribution of vector-valued Rademacher series,\sl\ Annals Prob.\
(to appear).\cr
F--R & A.~Fig\`a-Talamanca and D.~Rider,\rm\ A theorem of Littlewood and
lacunary series for compact groups,\sl\ Pacific J.\ Math.\ {\bf 16} (1966),
505--514.\cr
H--R & E.~Hewit and K.A.~Ross,\sl\ Abstract Harmonic Analysis II,\rm\ 
Springer--Verlag,
1970.\cr
Hi & P.~Hitczenko,\rm\ Domination inequality for martingale transforms of a 
Rademacher sequence,\sl\ Israel J.\ Math.\ (to appear).\cr
Ho & T.~Holmstedt,\rm\ Interpolation of quasi-normed spaces,\sl\ Math.\
Scand.\ {\bf 26} (1970), 177--199.\cr
Ka & J-P.~Kahane,\sl\ Some Random Series of Functions,\rm\ (2nd.\ Ed.) 
Cambridge studies
in advanced mathematics 5, 1985.\cr
Kh & A.~Khintchine,\rm\ \"Uber dyadische Br\"uche,\sl\ Math.\ Z.\ {\bf 18}
(1923), 109--116.\cr
Ko & A.N.~Kolmogorov,\rm\ \"Uber das Gesetz des iterieten Logarithmus,\sl\ 
Math.\ Ann.\ {\bf 101} (1929), 126--135.\cr
L--T & M.~Ledoux and M.~Talagrand,\sl\ Isoperimetry and Processes in
Probability in a Banach Space,\rm\ Springer--Verlag, 1991.\cr
M--P & M.B.~Marcus and G.~Pisier,\sl\ Random Fourier Series with Applications
to Harmonic Analysis,\rm\ Princeton University Press, 1981.\cr
Mo & S.J.~Montgomery-Smith,\rm\ The distribution of Rademacher sums,\sl\ Proc.\ 
A.M.S.\ {\bf 109}
(1990), 517--522.\cr
P1 & G.~Pisier,\rm\ De nouvelles caract\'erisations des ensembles de
Sidon,\rm\ Mathematical Analysis and Applications,\sl\ Advances in Math.\
Suppl.\ Stud., 7B (1981), 686--725.\cr
P2 & G.~Pisier,\sl\ Probabilistic Methods in the Geometry of Banach
Spaces,\rm\
Springer--Verlag, 1986.\cr
R--S & V.A.~Rodin and E.M.~Semyonov,\rm\ Rademacher series in symmetric
spaces,\rm\ Analyse Math.\ {\bf 1} (1975), 207--222.\cr
}

\bye